\documentclass{article}
\usepackage[usenames]{color}
\usepackage{amssymb}
\usepackage{amsmath}
\usepackage{amsthm}
\usepackage{amsfonts}
\usepackage{amscd}
\usepackage{graphicx}
\usepackage{diagbox}

\usepackage[colorlinks=true,
linkcolor=webgreen,
filecolor=webbrown,
citecolor=webgreen]{hyperref}

\definecolor{webgreen}{rgb}{0,.5,0}
\definecolor{webbrown}{rgb}{.6,0,0}

\usepackage[usenames]{color}
\usepackage{fullpage}
\usepackage{float}
\usepackage{graphics}
\usepackage{latexsym}
\usepackage{epsf}
\usepackage{breakurl}

\newcommand{\seqnum}[1]{\href{https://oeis.org/#1}{\rm \underline{#1}}}

\begin{document}

\theoremstyle{plain}
\newtheorem{theorem}{Theorem}
\newtheorem{proposition}[theorem]{Proposition}

\theoremstyle{definition}
\newtheorem{definition}[theorem]{Definition}
\newtheorem{example}[theorem]{Example}

\begin{center}
\vskip 1cm{\LARGE\bf Rational Dyck Paths\\

}

\vskip 1cm
\large
{Elena Barcucci, Antonio Bernini,
 Stefano Bilotta, Renzo Pinzani}\\
Dipartimento di Matematica e Informatica ``Ulisse Dini''\\
Universit\`a di Firenze\\

Viale G. B. Morgagni 65,
50134 Firenze, Italy\\

\end{center}

\vskip .2 in
\begin{abstract}
Given a positive rational $q$, we consider Dyck paths having height at most two with some constraints on the number of consecutive peaks and consecutive valleys, depending on $q$. We introduce a general class of Dyck paths, called rational Dyck paths, and provide the associated generating function, according to their semilength, as well as the construction of such a class. Moreover, we characterize some subsets of the rational Dyck paths that are enumerated by the $\mathbb Q$-bonacci numbers.
\end{abstract}

\section{Introduction}
In \cite{BKV} the authors introduced a new class of binary strings, Fibonacci $q$-decreasing words, which are enumerated by the Fibonacci generalized numbers. In \cite{BBP1}, later expanded in \cite{BBP2}, a new class of Dyck paths ({\it $q$-Dyck paths}, for short) enumerated by the Fibonacci generalized numbers was introduced. More recently, Kirgizov introduced a class of words, the $\mathbb Q$-bonacci words \cite{Ki}, larger than the one in \cite{BKV}, and enumerated by the $\mathbb Q$-bonacci numbers. In \cite{BBBP} we defined a class of Dyck paths ({\it $Q$-bonacci paths}, for short, from now on) larger than the one in \cite{BBP1,BBP2} and enumerated by the $\mathbb Q$-bonacci numbers.

In the present paper we define a new class of Dyck paths, {\it rational Dyck paths}. We give their construction and their generating function (which results to be rational). We also provide a very simple bijection between the the rational Dyck paths and the compositions having parts in a set \cite{HM}.

Moreover, since rational Dyck paths contain a subset of paths which are counted by the $\mathbb Q$-bonacci numbers, we give a bijection between them and a subset of corresponding paths in the set of $Q$-bonacci paths.

Finally, we show that there cannot exist a bijection between the rest of the rational Dyck paths and the rest of the $Q$-bonacci paths.

\section{Definition}
A Dyck path $P$ is a lattice path in $\mathbb {Z} ^{2}$ from 
$(0,0)$ to $(2n,0)$ with steps in $S=\lbrace (1,1),(1,-1)\rbrace$ that never passes below the $x$-axis. We indicate by $U$ the up step $(1,1)$ and by $D$ the down step $(1,-1)$, so that $P$ can be indicated by a suitable string over the alphabet $\{U,D\}.$

The height of $P$ is given by the maximum ordinate reached by its $U$ steps. In the paper we consider the set $\mathcal D$ of Dyck paths having heigth at most $2$.

The length of $P$ is the number of its steps and with $D_n$ we indicate the subset of $\mathcal D$ of the Dyck paths having semilength $n$ (i.e. the length divided by $2$). The empty path having semilength $0$ is denoted by $\varepsilon$.
Clearly, we have $\mathcal D=\displaystyle \bigcup_{n\geq 0} D_n$.

A \emph{1-peak} of $P$ is a substring $UD$ of $P$ where $U$ reaches ordinate $1$, and a \emph{0--valley} is a substring $DU$ of $P$ where $D$ touches the $x$-axis. 

If $L_1$ and $L_2$ are two subsets of $\mathcal D$, then the concatenation $L_1\cdot L_2$ is the set of all the paths of the form $PQ$ where $P\in L_1$ and $Q\in L_2$.

\section{Construction and generating function}\label{Section3}
A path $P\in \mathcal D$ $(P\neq \varepsilon)$ can be factorized 
highlighting the blocks of maximal length of consecutive $1$-peaks
and consecutive $0$-valleys:
\begin{equation}\label{factorization}
P=U(UD)^{p_1}(DU)^{v_1}(UD)^{p_2}(DU)^{v_2}\ldots(UD)^{p_k}(DU)^{v_k}D
\end{equation}
where possibly $p_1=0$ or $v_k=0$ (possibly both), where $(UD)^j$ (resp. $(DU)^j$) is the concatenation of $(UD)$ (resp. $(DU)$) with itself $j$ times. Note that if $k=0$, then $P=UD$.

Let $q=r/s$ be a positive rational $(q\in \mathbb Q^+)$ where the integers $r$ and $s$ are coprime.

\begin{definition}\label{Def1}
Let $ R_n^{r/s}$ be the set of Dyck paths $P\in\mathcal D$ of semilength $n$ where
$$
\frac{p_i}{v_i}\leq \frac{r}{s}
\text{ or equivalently $v_i\geq \left\lceil p_i\cdot\frac{s}{r}\right\rceil$}
$$
for $i=1,2,\ldots,k.$

Moreover, let $\mathcal R^{r/s}$ be the class of the \emph{rational Dyck paths} collecting the sets $R_n^{r/s}$, for each $n$:
$$
\mathcal R^{r/s}=\bigcup_{n\geq 0}R_n^{r/s}.
$$ 
\end{definition}
In other words, the rational Dyck paths of $\mathcal R^{r/s}$ are such that if $p$ consecutive $1$-peaks occur, then they are immediately followed by at least
$\left\lceil p\cdot\frac{s}{r}\right\rceil$ consecutive $0$-valleys.  
Note that if $P\in\mathcal R^{r/s}$, then in its factorization (\ref{factorization}) it can not be $v_k=0$.

\begin{example}
If $q=3/4$ and
$$
P_1=UUDDUDUUDUDDUDUD=U(UD)(DU)^2(UD)^2(DU)^2D,
$$
then $k=2, p_1=1, v_1=2, p_2=2, v_2=2$. According to Definition \ref{Def1}, it should be $v_2\geq3$. Therefore $P_1\notin \mathcal R^{3/4}$.

The path
$$
P_2=UUDDUDUUDUDUDUDD=U(UD)(DU)^2(UD)^4D
$$
is not allowed $(P_2\notin \mathcal R^{3/4})$ since in this case $v_2=0$ (there are no valleys after the last block of $1$-peaks).

The path
$$
P_3=UUDDUDUUDUDDUDUDUDUD=U(UD)(DU)^2(UD)^2(DU)^4D
$$
is allowed.
\end{example}
Clearly, a path $P$ different from the empty set $\varepsilon$ $\bigl(P\in \{\mathcal R^{r/s}\setminus \varepsilon\}\bigr)$ can be recursively constructed as follows:

\begin{itemize}
\item $P=UD\cdot P', \text{ with } P'\in\mathcal R^{r/s} \text{ (possibly $P'=\varepsilon$)}$, either
\item $P=U(UD)^p(DU)^{\lceil ps/r \rceil-1}D\cdot P', \text{ with } P'\in\{\mathcal R^{r/s}\setminus\varepsilon\}, \text{ and } p\geq 1$ .
\end{itemize}
Note that in the second case $P'\neq \varepsilon$, otherwise in the path $P$ the $p$ consecutive $1$-peaks of the prefix $U(UD)^p(DU)^{\lceil ps/r \rceil-1}D$ would not be followed by the right number of consecutive $0$-valleys which must be at least $\lceil ps/r \rceil$, according to Definition \ref{Def1}: the missing one (observe that $\lceil ps/r \rceil-1$ consecutive $0$-valleys occur in the prefix itself) is given by the final $D$ step of the prefix and the first $U$ step of $P'$.

Therefore, the construction of the rational Dyck paths $\mathcal R^{r/s}$ can be expressed by:
\begin{equation}\label{costruzione}
\mathcal R^{r/s}=\varepsilon\ \cup\ UD\cdot\mathcal R^{r/s}\ 
\bigcup_{p\geq 1} U(UD)^p(DU)^{\lceil ps/r \rceil-1}D\cdot\{\mathcal R^{r/s}\setminus\varepsilon\}
\end{equation}
We denote by $\delta ^{r/s}(x)$ the generating function of $\mathcal R^{r/s}$, according to the semilength of the paths:
$$
\delta ^{r/s}(x)=\sum_{P\in\mathcal R^{r/s}}x^{|P|}
$$
where $|P|$ denotes the semilength of $P$. From (\ref{costruzione}) it is possible to deduce the functional equation for the generating function $\delta ^{r/s}(x)$:
\begin{equation}\label{eq_funzionale}
\delta ^{r/s}(x)=1 + x\delta ^{r/s}(x) + 
\sum_{p\geq 1} x^{p+\lceil ps/r \rceil}\bigl(\delta ^{r/s}(x)-1\bigr).
\end{equation}
In the above formula the sum $\sum_{p\geq 1} x^{p+\lceil ps/r \rceil}$
can be manipulated so that the functional equation boils down to:
\begin{align}\label{funzione_generatrice}
\delta ^{r/s}(x)&=1 + x\delta^{r/s}(x) + 
\bigl(\delta^{r/s}(x)-1\bigr)
\sum_{k= 1}^{r} x^{k+\lceil ks/r \rceil}
\sum_{j\geq 0}x^{(r+s)j}\nonumber\\
&=1+ x\delta^{r/s}(x) + 
\bigl(\delta^{r/s}(x)-1\bigr)\cdot\frac{\displaystyle\sum_{k= 1}^{r} x^{k+\lceil ks/r \rceil}}{1-x^{r+s}}
\end{align}
and the generating function is:
\begin{align}\label{GEN_FUN}
\delta^{r/s}(x)&=\frac{1-x^{r+s}-\displaystyle\sum_{k= 1}^{r} x^{k+\lceil ks/r \rceil}}{(1-x)(1-x^{r+s})-\displaystyle\sum_{k= 1}^{r} x^{k+\lceil ks/r \rceil}}\nonumber\\
&=\frac{1-\displaystyle\sum_{k= 1}^{r-1} x^{k+\lceil ks/r \rceil}-2x^{r+s}}{1-x-\displaystyle\sum_{k= 1}^{r-1} x^{k+\lceil ks/r \rceil}-2x^{r+s}+x^{r+s+1}}.
\end{align}
Following the above construction summarized in (\ref{costruzione}) it is not difficult to realize that the generation of the sets
$R_n^{r/s}$ is given by:
\begin{equation}
R_n^{r/s}=
\begin{cases}
\varepsilon,  & \text{if $n=0$};\\
&\\
UD\cdot R_{n-1}^{r/s}\displaystyle\bigcup_{p\geq1} \{Y_{n,p}\}, & \text{if $n\geq1$};\\
\end{cases} 
\end{equation}
where
\begin{equation*}
Y_{n,p}=
\begin{cases}
U(UD)^{p}(DU)^{\lceil{p s/r}\rceil-1}D\cdot R_{n-p-\lceil{p s/r}\rceil}  & \text{if $n-p-\lceil{p s/r}\rceil\geq1$};\\
&\\
\emptyset \text{ (empty set)},  & \text{otherwise}.\\
\end{cases}
\end{equation*}
Hence, denoting by $w_n^{r/s}$ the cardinality of
$R_n^{r/s}$,  the recurrence relation for $w_n^{r/s}$ is:
\begin{equation}
w_n^{r/s}=
\begin{cases}
1,  & \text{if $n=0$};\\
&\\
w_{n-1}^{r/s}+ \displaystyle\sum_{p\geq 1}u_{n,p}, & \text{if $n\geq1$};\\
\end{cases} 
\end{equation}
where
\begin{equation*}
u_{n,p}=
\begin{cases}
w_{n-p-\lceil{p s/r}\rceil}  & \text{if $n-p-\lceil{p s/r}\rceil\geq1$};\\
&\\
0,  & \text{otherwise}.\\
\end{cases}
\end{equation*}

\section{Links with other structures}\label{LINK}
If $r=s=1$ the generating function (\ref{GEN_FUN}) becomes
$$
\delta (x)=\frac{1-2x^2}{1-x-2x^2+x^3}
$$
corresponding to the sequence $1,1,1,2,3,6,10,19,33,\ldots$ which is a shift of the \seqnum{A028495} sequence in The On-line Encyclopedia of Integer Sequences \cite{S}. The sequence enumerates the compositions of $n$ with parts into $A=\{1,2,4,6,\ldots\}$ according to their length,  as we can deduce from \cite{HM} and \cite{BRV}.

On closer inspection, our sets $R_{n+1}^{r/s}$ are in bijection with the sets of compositions $C_n^{A_{r/s}}$ of $n$ with parts into the set
$$A_{r/s}=\{1\}\cup\{p+\lceil ps/r\rceil \mid p\geq 1\}.$$
The bijection is accomplished simply mapping the part $1$ with the prefix $UD$, and the part $p+\lceil ps/r\rceil$ with prefix $U(UD)^{p}(DU)^{\lceil{p s/r}\rceil-1}D$. Finally, adding the factor $UD$ at the end of the obtained Dyck path.
For example, if $q=3/4$ then $A_{3/4}=\{{1,3,5,7,10,12,14,\ldots}\}$ and the correspondence between $C_n^{A_{3/4}}$ and $R_{n+1}^{3/4}$ can be seen in Table \ref{Table1}.

\begin{table}[ht]
\centering
\begin{tabular}{ | l | c | c |}
       \hline 
          & $C_n^{A_{3/4}}$ & $R_{n+1}^{3/4}$ \\ \hline
    $n=0$ & $\varepsilon$ & $UD$ \\ \hline
    $n=1$ & $1$ & $UDUD$ \\ \hline
    $n=2$ & $11$ & $UDUDUD$ \\ \hline
    $n=3$ & $111,\ 3$ & $UDUDUDUD$, $UUDDUDUD$ \\ \hline
 		  & $1111$ & $UDUDUDUDUD$\\	   
    $n=4$ & $13$ & $UDUUDDUDUD$ \\ 
          & $31$ & $UUDDUDUDUD$ \\ 
    \hline
    	  & $11111$ & $UDUDUDUDUDUD$ \\    
    	  & $113$   & $UDUDUUDDUDUD$ \\	   
    $n=5$ & $131$   & $UDUUDDUDUDUD$ \\ 
          & $311$   & $UUDDUDUDUDUD$ \\ 
    	  & $5$     & $UUDUDDUDUDUD$ \\
    \hline
  $\cdots$& $\cdots$& $\cdots$\\   
    \hline    
\end{tabular}
\caption{$C_n^{A_{3/4}}$ and $R_{n+1}^{3/4}$}
\label{Table1}
\end{table}
In a recent paper by Kirgizov \cite{Ki} the $q$-decreasing words, previously introduced with $q\in \mathbb N^+$  \cite{BKV}, have been generalized to the case where $q\in\mathbb Q^+$.
If $q\in \mathbb N^+$, the $q$-decreasing words are enumerated \cite{BKV} by the well known $q$-generalized Fibonacci numbers \cite{F,M}, while in the case where $q\in \mathbb Q^+$, they are enumerated by the so called $\mathbb Q$-bonacci numbers \cite{Ki}.

For some particular values $q=r/s$, the $q$-decreasing words of length $n$ are enumerated by some sequences that, when properly filled with 1's at the beginning, coincide with the ones enumerating some restricted compositions of $n$ with a finite number of parts depending on $r/s$.
For example the A060691, A117760 sequences in The On-line Encyclopedia of Integer Sequences \cite{S} corresponding to the values $q=2/3$, and $q=3/4$ give the enumeration of the $q$-decreasing words and the enumeration of the compositions of $n$ with parts in $\{1,3,5\}$ and $\{1,3,5,7\}$, respectively.

More precisely, the involved compositions are the ones where the summands of $n$ belong to the finite set 
$$
\widetilde{A}_{r/s}=\{1\}\cup\{p+\lceil ps/r\rceil \mid 1\leq p\leq r\}.
$$

This particular background where restricted compositions, Dyck paths, and $q$-decreasing words seem to be connected in an intriguing way, prompted us to investigate on the set $\mathcal R^{r/s}$  in order to  
\begin{itemize}
\item find Dyck paths in bijection with $C_n^{\widetilde{A}_{r/s}}.$
\item find Dyck paths having the same enumeration of the $q$-decreasing words according to the semilength, for $q\in\mathbb Q^+$;
\end{itemize}

\subsection{Dyck paths and restricted compositions}\label{COMP}
Referring again to the factorization (\ref{factorization}) 
\begin{equation*}
	P=U(UD)^{p_1}(DU)^{v_1}(UD)^{p_2}(DU)^{v_2}\ldots(UD)^{p_k}(DU)^{v_k}D
\end{equation*}
for Dyck paths $P\in\mathcal D$ we insert a constraint on the number 
consecutive $1$-peaks in each block. 

\begin{definition}\label{Def3}
	Let $\widetilde R_n^{r/s}$ be the set of Dyck paths $P\in\mathcal D$ of semilength $n$ where:
	\begin{itemize}
		\item $p_i\leq r$ for $i=1,2,\ldots,k$;
		\item for $i=1,2,\ldots,k$ it has to be
		$$
		\frac{p_i}{v_i}\leq \frac{r}{s}
		\text{ or equivalently $v_i\geq \left\lceil p_i\cdot\frac{s}{r}\right\rceil.$}
		$$
	\end{itemize}
	Moreover, let $\widetilde R^{r/s}$ be the class collecting the sets $\widetilde R_n^{r/s}$, for each $n$:
	$$
	\widetilde{\mathcal R}^{r/s}=\bigcup_{n\geq 0}\widetilde R_n^{r/s}.
	$$ 
\end{definition}
All the steps showed for $\mathcal R^{r/s}$ in Section \ref{Section3} can be easily adjusted to Definition \ref{Def3}. More precisely:
\begin{itemize}
	\item the construction of $\widetilde{\mathcal R}^{r/s}$ can be described by
	\begin{equation*}
		\widetilde{\mathcal R}^{r/s}=\varepsilon\ \cup\ UD\cdot
		\widetilde{\mathcal R}^{r/s}\ 
		\bigcup_{p=1}^r U(UD)^p(DU)^{\lceil ps/r \rceil-1}D\cdot
		\{\widetilde{\mathcal R}^{r/s}\setminus\varepsilon\};
	\end{equation*}
	\item the corresponding functional equation for the generating function $\widetilde{\delta}^{r/s}(x)$ of $\widetilde{\mathcal R}^{r/s}$ is
	\begin{equation*}
		\widetilde{\delta}^{r/s}(x)=1 + x\widetilde{\delta}^{r/s}(x) + 
		\sum_{p=1}^{r} x^{p+\lceil ps/r \rceil}\bigl(\widetilde{\delta}^{r/s}(x)-1\bigr);
	\end{equation*}
	\item the generating function is
	\begin{equation}\label{gen_fun_tilde}
		\widetilde{\delta}^{r/s}(x)=\frac{1-\displaystyle\sum_{p=1}^{r} x^{p+\lceil ps/r \rceil}}{1-x-\displaystyle\sum_{p= 1}^{r} x^{p+\lceil ps/r \rceil}};
	\end{equation}
	\item the generation of the sets $\widetilde R_n^{r/s}$ is given by:
	\begin{equation*}
		\widetilde R_n^{r/s}=
		\begin{cases}
			\varepsilon,  & \text{if $n=0$};\\
			&\\
			UD\cdot\widetilde R_{n-1}^{r/s}\displaystyle\bigcup_{p=1}^r \{\widetilde{Y}_{n,p}\}, & \text{if $n\geq1$};\\
		\end{cases} 
	\end{equation*}
	where
	\begin{equation*}
		\widetilde{Y}_{n,p}=
		\begin{cases}
			U(UD)^{p}(DU)^{\lceil{p s/r}\rceil-1}D\cdot\widetilde R_{n-p-\lceil{p s/r}\rceil}  & \text{if $n-p-\lceil{p s/r}\rceil\geq1$};\\
			&\\
			\emptyset \text{ (empty set)},  & \text{otherwise}.\\
		\end{cases}
	\end{equation*}
	Hence, denoting by $\widetilde{w}_n^{r/s}$ the cardinality of
	$\widetilde R_n^{r/s}$,  the recurrence relation for $\widetilde{w}_n^{r/s}$ is:
	\begin{equation}\label{rec_rel_tilde}
		\widetilde{w}_n^{r/s}=
		\begin{cases}
			1,  & \text{if $n=0$};\\
			&\\
			\widetilde{w}_{n-1}^{r/s}+ \displaystyle\sum_{p=1}^r\widetilde{u}_{n,p}, & \text{if $n\geq1$};\\
		\end{cases} 
	\end{equation}
	where
	\begin{equation*}
		\widetilde{u}_{n,p}=
		\begin{cases}
			\widetilde{w}^{r/s}_{n-p-\lceil{p s/r}\rceil}  & \text{if $n-p-\lceil{p s/r}\rceil\geq1$};\\
			&\\
			0,  & \text{otherwise}.\\
		\end{cases}
	\end{equation*}
\end{itemize}
The bijection between the sets $C_n^{\widetilde{A}_{r/s}}$
and $\widetilde R_{n+1}^{r/s}$ is equal to the one between the sets $C_n^{A_{r/s}}$ and $R_{n+1}^{r/s}$, defined at the beginning of Section \ref{LINK}. Here, we only note that the set $\widetilde{A}_{r/s}$ containing the parts for the composition of $n$ is a finite set, corresponding to the fact that the paths $P\in\widetilde R_n^{r/s}$ can not have more than $r$ consecutive $1$-peaks. Moreover, thanks to the bijection, we observe that $\widetilde{w}_n^{r/s}$ is the number of the composition of $n-1$ for $n\geq 1$ with parts into $\widetilde{A}^{r/s}$. This is the reason why the sequence defined by (\ref{rec_rel_tilde}) does not match exactly
the one enumerating the composition of $n$ with parts into $\widetilde{A}^{r/s}$. Specifically, the sequence $(\widetilde{w}_n^{r/s})_{n\geq 0}$ is obtained by the one enumerating the compositions by inserting a $1$ at the beginning. Indeed, the generating function
$$
\frac{1}{1-\displaystyle\sum_{\ell\in\widetilde{A}^{r/s}}x^{\ell}}
$$
for the compositions of $n$ with parts in $\widetilde{A}^{r/s}$ \cite{HM} can be obtained as $\widetilde{\delta}^{r/s}(x)-1$.

\subsection{Dyck paths and \texorpdfstring{$q$}{q}-decreasing words}
Here we deal with the family of Dyck paths having the same enumeration, according to their semilength, of the $q$-decreasing words as generalized by Kirgizov \cite{Ki}.
These paths are the ones introduced in \cite{BBBP} where they were presented with a different but equivalent approach. Here, after their definition, following a similar argument used in Section \ref{Section3} we recall only the fundamental points in order to arrive to the enumerating sequences.

Recalling the factorization (\ref{factorization}) 
\begin{equation*}
P=U(UD)^{p_1}(DU)^{v_1}(UD)^{p_2}(DU)^{v_2}\ldots(UD)^{p_k}(DU)^{v_k}D
\end{equation*}
for Dyck paths $P\in\mathcal D$ we keep the restriction on the number of consecutive $1$-peaks and relax the assumption on $v_k$.

\begin{definition}\label{decreasing}
Let ${\text Q}_n^{r/s}$ be the set of Dyck paths $P\in\mathcal D$ of semilength $n$ where:

\begin{itemize}
\item $p_i\leq r$ for $i=1,2,\ldots,k$;
\item for $i=1,2,\ldots,k-1$ it has to be
$$
\frac{p_i}{v_i}\leq \frac{r}{s}
\text{\qquad  or equivalently $v_i\geq \left\lceil p_i\cdot\frac{s}{r}\right\rceil;$}
$$
\item for $i=k$
\begin{itemize}
\item if $p_k=r$, then $v_k\geq 0$;
\item if $p_k<r$, then
 either $v_k=0$ or $v_k\geq \left\lceil p_k\cdot\frac{s}{r}\right\rceil.$
\end{itemize}
\end{itemize}
Moreover, let $\mathcal Q^{r/s}$ be the class of \emph{$Q$-bonacci paths} collecting the sets $Q_n^{r/s}$, for each $n$:
$$
\mathcal Q^{r/s}=\bigcup_{n\geq 0}Q_n^{r/s}.
$$ 
\end{definition}
In other  words, in a path $P\in Q_n^{r/s}$ a block $B$ of $p$ consecutive $1$-peaks must be followed by at least $\lceil p\cdot\frac{s}{r}\rceil$ consecutive $0$-valleys. The only case where a block $B$ can be followed by less than $\lceil p\cdot\frac{s}{r}\rceil$ consecutive $0$-valleys is when $B$ is the last one and $B=(UD)^r$. Moreover, the $p$ consecutive $1$-peaks in each block $B$ cannot be more than 
$r$.

\begin{example}
If $q=4/5$ and
$$
P_1=UUDDUDUUDUDDUDUD=U(UD)(DU)^2(UD)^2(DU)^2D,
$$
then $k=2, p_1=1, v_1=2, p_2=2, v_2=2$. According to Definition \ref{decreasing}, it should be $v_2\geq3$. Therefore $P_1\notin \mathcal Q^{4/5}$.

The path
$$
P_2=UUDDUDUUDUDUDUDDUDUD=U(UD)(DU)^2(UD)^4(DU)^2D
$$
is allowed $(P_2\in \mathcal Q^{4/5})$ since in this case $p_2=4=r$, so that it is only required $v_2\geq 0$.

The path
$$
P_3=UUDDUDUUDUDUDD=U(UD)(DU)^2(UD)^3D
$$
is allowed since in this case $p_2=3<r$ and $v_2= 0$.
\end{example}
These paths are the ones introduced in \cite{BBBP} where they were presented with a different but equivalent definition. Here we only recall the fundamental points in order to arrive to the enumerating sequences.

\begin{itemize}
	\item The construction of $\mathcal Q^{r/s}$ can be summarized by
	\begin{align*}
		\mathcal Q^{r/s}=&\varepsilon\ \cup\ UD\cdot
		\mathcal Q^{r/s}\ 
		\bigcup_{\ell=1}^{r+s-1} Upr_{\ell}\bigl((UD)^r(DU)^{s-1}\bigr)D\\
		&\bigcup_{p=1}^r U(UD)^p(DU)^{\lceil ps/r \rceil-1}D\cdot
		\{\mathcal Q^{r/s}\setminus\varepsilon\}, 
	\end{align*}
	where, if $A$ is a path, then $pr_{\ell}(A)$ denotes the prefix of semilength $\ell$ of $A$.
	
	Specifically, a path $P\in\mathcal Q^{r/s}$, as in the case where $P\in \widetilde{\mathcal R}^{r/s}$ (Section \ref{COMP}),  starts with the factor $UD$, or one of the factors $U(UD)^p(DU)^{\lceil ps/r \rceil-1}D$ with $p=1,2,\ldots,r$, concatenated with a suitable path of $\mathcal Q^{r/s}$. Here ($P\in\mathcal Q^{r/s}$) the path $P$ can be also a prefix of $(UD)^p(DU)^{\lceil ps/r \rceil-1}$ preceded by a step $U$ and followed by a step $D$.     
	\item The functional equation for the generating function $\chi^{r/s}(x)$ of $\mathcal Q^{r/s}$ is:
	\begin{equation*}
		\chi^{r/s}(x)=1 + x\chi^{r/s}(x) +
		\sum_{j=2}^{r+s}x^j\  
		+\sum_{p=1}^{r} x^{p+\lceil ps/r \rceil}\bigl(\chi^{r/s}(x)-1\bigr),
	\end{equation*}
where the last sum takes track of the paths given by the prefixes of 
$(UD)^r(DU)^{s-1}$ preceded by $U$ and followed by $D$.	
	\item The generating function is
	\begin{equation*}
		\chi^{r/s}(x)=\frac
		{1+\displaystyle\sum_{j=2}^{r+s}x^j-\sum_{p=1}^{r}x^{p+\lceil ps/r \rceil}}
		{1-x-\displaystyle\sum_{p=1}^{r}x^{p+\lceil ps/r \rceil}}\ .
	\end{equation*}
	\item The generation of the sets
	$Q_n^{r/s}$ is described by:
	
	\begin{equation*}
		Q_n^{r/s}=
		\begin{cases}
			\varepsilon,  & \text{if $n=0$};\\
			&\\
			UD,&\text{if $n=1$};\\
			&\\
			UD\cdot Q_{n-1}^{r/s}
			\cup U pr_{n-1}\bigl((UD)^r(DU)^{s-1}\bigr)D		
			\displaystyle\bigcup_{p=1}^r \{S_{n,p}\}, & \text{if $n\geq2$};\\
		\end{cases} 
	\end{equation*}
	where
	\begin{equation*}
		S_{n,p}=
		\begin{cases}
			U(UD)^{p}(DU)^{\lceil{p s/r}\rceil-1}D\cdot Q_{n-p-\lceil{p s/r}\rceil}  & \text{if $n-p-\lceil{p s/r}\rceil\geq1$};\\
			&\\
			\emptyset \text{ (empty set)},  & \text{otherwise}.\\
		\end{cases}
	\end{equation*}
	Hence, denoting by $v_n^{r/s}$ the cardinality of
	$Q_n^{r/s}$,  the recurrence relation for $v_n^{r/s}$ is:
	\begin{equation}
		v_n^{r/s}=
		\begin{cases}
			1,  & \text{if $n=0$};\\
			&\\
			v_{n-1}^{r/s}+ \displaystyle\sum_{p=1}^rt_{n,p}, & \text{if $n\geq1$};\\
		\end{cases} 
	\end{equation}
	where
	\begin{equation*}
		t_{n,p}=
		\begin{cases}
			v^{r/s}_{n-p-\lceil{p s/r}\rceil}  & \text{if $n-p-\lceil{p s/r}\rceil\geq1$};\\
			&\\
			0,  & \text{otherwise}.\\
		\end{cases}
	\end{equation*}
\end{itemize} 
Denoting by $W_{r/s}(x)$ the generating function for the
$q$-decreasing words as defined in \cite{Ki}, it is routine to check that
$$
\chi^{r/s}(x)=1+xW_{r/s}(x)\ ,
$$
so that the sequences $\bigl(v_n^{r/s}\bigr)_{n\geq 0}$ are obtained by the ones enumerating the $q$-decreasing words by inserting a $1$ at the beginning. Hence, all the sequences obtained in \cite{Ki} have a different combinatorial interpretation by means of $Q$-bonacci paths.

\subsection{The classes \texorpdfstring{$\mathcal Q^{r/s}$}{Qrs} and \texorpdfstring{$\widetilde{\mathcal R}^{r/s}$}{Drs}
}
In this section, we are going to prove that for particular cases, some subsets of $\mathcal Q^{r/s}$ and $\widetilde{\mathcal R}^{r/s}$ are in biejection. More precisely, we have the following proposition.

\begin{proposition}
	If $s=tr+1$ with $t \in \mathbb{N}$, then there exists a biejction $\phi: \emph{{\text Q}}_n^{r/s} \rightarrow \widetilde{R}_{n+t+1}^{r/s}$, for $n \geq 0$.
\end{proposition}   

\begin{proof}	
	Let be $P \in Q_n^{r/s}$. The general idea in order to obtain the corresponding path $\phi(P)$ in $\widetilde R_{n+t+1}^{r/s}$ is adding the suffix $(UD)^{t+1}$ (so that $\phi(P)$ ends with at least $t+1$ consecutive $0$-valleys and its semilength is exactly $n+t+1$) and eventually replacing with $0$-valleys a suitable number of the $p_k$ rightmost consecutive $1$-peaks of $P$.
	
	Referring to factorization (\ref{factorization}), if $p_k = j$, let $\nu_j=\left\lceil j\cdot\frac{s}{r}\right\rceil$ so that $v_k \geq \nu_j$, according to Definition \ref{Def3}.
	
	It is easy to check that 
	\begin{equation}\label{nu}
		\nu_j=jt+1,
	\end{equation}
	being $1 \leq j \leq r$.
	
	We are going to consider the following cases:
	\begin{itemize}
		\item If $P \in \widetilde{R}_{n}^{r/s} (\subseteq Q_n^{r/s})$, then it is easily seen that $\phi(P)=P (UD)^{t+1} \in \widetilde{R}_{n+t+1}^{r/s}$;
		\item If $P=\alpha U (UD)^j D \in Q_n^{r/s}$, where $\alpha \in \widetilde{R}_{n-j-1}^{r/s}$, then the path $\phi(P)$ can be expressed by $\phi(P)=\alpha U (UD)^{j-h} (DU)^h D (UD)^{t+1}$ with $h$ is the minimum integer such that $h+t+1 \geq \nu_{j-h}= t(j-h)+1$, taking into account formula (\ref{nu}). Hence $h =\left\lceil \frac{jt+1}{t+1} -1\right\rceil $. 
		
		Suppose that there exists $P' \in \widetilde{R}_n^{r/s}$ such that $\phi(P)=\phi(P')$, then we can write $P'=\alpha U (UD)^{j-h} (DU)^h D$. According to formula (\ref{nu}), we have $h \geq (j-h)t+1$. So that, it should be $h \geq \left\lceil \frac{jt+1}{t+1}\right\rceil > \left\lceil \frac{jt+1}{t+1} -1\right\rceil$, against the preceding value of $h$.
		
		\item If $P=\alpha U (UD)^r (DU)^j D$, with $1 \leq j \leq s-1$, then we distinguish 2 cases: 
		\begin{itemize}
			\item If $j \geq s-t-1$, then $\phi(P)=P (UD)^{t+1} \in \widetilde{R}_{n+t+1}^{r/s}$.
			\item If $1 \leq j < s-t-1$, then $\phi(P)=\alpha U (UD)^{r-h} (DU)^{j+h} D (UD)^{t+1}$, with $h$ is the minimum integer such that $j+h+t+1 \geq \nu_{r-h}= (r-h)t+1$ thanks to (\ref{nu}), hence $h =\left\lceil \frac{rt-j+1}{t+1} -1\right\rceil $.    
		\end{itemize}
	With similar arguments described in the previous bullet, one can show that there not exists $P' \in \widetilde{R}_n^{r/s}$ such that $\phi(P)=\phi(P')$.  
	\end{itemize}	
	
	The inverse map $\phi^{-1}$ can be obtained with similar considerations. Given $P \in \widetilde{R}_{n+t+1}^{r/s}$ of the form $P=\alpha U (UD)^j (DU)^h D$, with $h \geq \nu_j \geq t+1 $ and $\alpha \in \widetilde{R}_{n+t-j-h}^{r/s}$, we are going to consider the following cases:
	
	\begin{itemize}
		\item If $h-(t+1) \geq \nu_j$, then $\phi^{-1}(P)=\alpha U (UD)^j (DU)^{h-t-1} D$, which belongs to $\widetilde{R}_{n}^{r/s} (\subseteq Q_n^{r/s})$.
		\item If $h-(t+1) < \nu_j$ and $j+h-t-1 \leq r$, then $\phi^{-1}(P)=\alpha U (UD)^{j+h-t-1}D$, which belongs to $Q_n^{r/s}$.
		\item If $h-(t+1) < \nu_j$ and $j+h-t-1 > r$, then $\phi^{-1}(P)=\alpha U (UD)^r (DU)^{j+h-t-1-r} D$, which belongs to $Q_n^{r/s}$.
	\end{itemize}    
	It is routine to show that $\phi(\phi^{-1}(P))=P$. The map $\phi$ is the required bijection. Note that, the hypothesis $s=tr+1$ is not used for the definition of $\phi^{-1}$.
\end{proof}
If $s \neq tr + 1$, then $tr+1 < s < (t+1)r$. In this case, we have $t\leq \nu_{j+1} - \nu_j \leq t+1$ and of course there exists $j$ such that $\nu_{j+1} + \nu_j=t+1$. Let $j_0$ be the minimum integer such that $\nu_{r-j_0} -\nu_{r-(j_0+1)}=t+1$. We have $\nu_{r-j_0} = s-tj_0$ and $\nu_{r-(j_0+1)}=s-t(j_0+1)-1$. We how show that there exist two paths $P,P'\in \mathcal Q_n^{r/s}$ such that $\phi(P)=\phi(P')$. Let
$$
P=\alpha U (UD)^r (DU)^{s-(t+1)(j_0+1)-1}D \in  Q_n^{r/s} \setminus \widetilde{R}_{n}^{r/s}
$$
and
$$
P'=\alpha U (UD)^{r-(j_0+1)} (DU)^{s-t(j_0+1)-1}D \in \widetilde{R}_{n}^{r/s},
$$
by applying the map $\phi$, we obtain $\phi(P)=\phi(P')=\alpha U (UD)^{r-(j_0+1)} (DU)^{s-t(j_0+1)+t}D$. Hence, the map $\phi$ is not injective but it is surjective. Therefore, the sets $\widetilde{R}_{n+t+1}^{r/s}$ and $Q_n^{r/s}$ have different cardinalities and the following proposition is proved:

\begin{proposition}
	The sets $\widetilde{R}_{n+t+1}^{r/s}$ and $Q_n^{r/s}$ are in biejection if and only if $s=tr+1$.
\end{proposition}
We point out that in the the class of $Q$-bonacci paths there are some rational Dyck paths. They are exactly the paths of the class $\widetilde{\mathcal R}^{r/s}$. In other words

$$
\mathcal R^{r/s}\cap\mathcal Q^{r/s} =\widetilde{\mathcal R}^{r/s}\ .
$$

\section{Conclusions}
In the present work we consider the construction and the enumeration of the class $\mathcal R^{r/s}$, that is, the class of rational Dyck paths. We prove that some subsets of $\mathcal R^{r/s}$ are in biejection with subsets of the $Q$-bonacci paths (enumerated by $\mathbb Q$-bonacci numbers), if and only if a given relation between $r$ and $s$ is satisfied.

In our opinion there are some open questions that deserve to be answered. One of the most intriguing is the fact that $Q$-bonacci paths in \cite{BBBP} and $\mathbb{Q}$-bonacci words in \cite{Ki} are both enumerated by $\mathbb Q$-bonacci numbers, but an explicit bijection between the mentioned classes is not yet defined. 

A second one could deal with the definition of a Gray code for rational Dyck paths or/and $Q$-bonacci paths. 

%

%
%
%
%
%
%
%
%


\begin{thebibliography}{99}
\bibitem{BBBP} E. Barcucci, A, Bernini, S, Bilotta, and R. Pinzani, Dyck paths enumerated by the $\mathbb Q$-bonacci numbers,
in S. Brlek and L. Ferrari, Editors, {\it Electronic Proceedings in Theoretical Computer Science} {\bf 403} (2024), 49--53.

\bibitem{BBP1} E. Barcucci, A, Bernini, and R. Pinzani,
From the {F}ibonacci to {P}ell numbers and beyond via {D}yck paths,
{\it Pure Math. Appl. (PU.M.A.)}
{\bf 30}
(2022),	17--22.

\bibitem{BBP2} E. Barcucci, A, Bernini, and R. Pinzani,
Sequences from Fibonacci to Catalan: A combinatorial interpretation via Dyck paths,
{\it RAIRO Theor. Inform. Appl. (RAIRO:ITA)}
{\bf 58}
(2024),	\#8.


\bibitem{BKV} J.-L. Baril, S. Kirgizov, and V. Vajnovszki,
Gray codes for Fibonacci $q$-decreasing words,
{\it Theoret. Comput. Sci.} {\bf 297} (2022), 120--132.

\bibitem{BRV} J.-L. Baril, J. L. Ram\'irez, and F. A. Velandia, Bijections between directed-column convex polyominoes and
restricted compositions, {\it Theoret. Comput. Sci.} {\bf 1003} (2024), 114626.

\bibitem{F} M. Feinberg, Fibonacci-Tribonacci, {\it Fibonacci Quart.} {\bf 1} (1963), 71--74.


\bibitem{HM} S. Heubach, and T. Mansour, Compositions of $n$ with parts in a set, {\it Congressus Numerantium} {\bf 168} (2004), 127--143.

\bibitem{Ki} S. Kirgizov, $\mathbb{Q}$-bonacci words and numbers, {\it Fibonacci Quart.} {\bf 60} (2022), 187--195.

\bibitem {M} E. P. Miles Jr., (1960): Generalized Fibonacci numbers and associated matrices. {\it Amer. Math. Monthly}, {\bf 67} (1960), 745--752.

\bibitem{S} N. J. A. Sloane, The On-Line Encyclopedia of Integer Sequences, \url{https://oeis.org}.


\end{thebibliography}
\end{document}